\documentclass[10pt]{article}
\usepackage{mathrsfs}
\usepackage{amssymb,amsfonts,amsmath, psfrag,eepic,colordvi,epsfig}
\topskip 
\numberwithin{equation}{section}
\newtheorem{theorem}{Theorem}[section]

\newtheorem{conjecture}[theorem]{Conjecture}

\newtheorem{lemma}[theorem]{Lemma}

\begin{document}
\parskip 7pt

\pagenumbering{arabic}
\def\sof{\hfill\rule{2mm}{2mm}}
\def\ls{\leq}
\def\gs{\geq}
\def\SS{\mathcal S}
\def\qq{{\bold q}}
\def\MM{\mathcal M}
\def\TT{\mathcal T}
\def\EE{\mathcal E}
\def\lsp{\mbox{lsp}}
\def\rsp{\mbox{rsp}}
\def\pf{\noindent {\it Proof.} }
\def\mp{\mbox{pyramid}}
\def\mb{\mbox{block}}
\def\mc{\mbox{cross}}
\def\qed{\hfill \rule{4pt}{7pt}}
\def\pf{\noindent {\it Proof.} }
\textheight=22cm

\begin{center}
{\Large\bf Proofs of some conjectures
 of Chan-Mao-Osburn on
 Beck's partition  statistics}
\end{center}

\begin{center}

$^{1}$Liuxin Jin, $^{2}$Eric H. Liu and $^{3}$Ernest X.W. Xia

$^{1}$Department of Mathematics, \\
 Jiangsu University, \\
 Jiangsu, Zhenjiang, 212013, P. R. China\\[6pt]

$^{2}$School of Statistics and Information, \\
 Shanghai University of
International Business and Economics,\\
  Shanghai, 201620, P. R. China\\[6pt]

$^3$School of Mathematical Sciences, \\
  Suzhou University of Science and
Technology, \\
 Suzhou,  215009, Jiangsu Province,
 P. R. China

Email: liuxj@ujs.edu.cn, liuhai@suibe.edu.cn,
  ernestxwxia@163.com

\end{center}

%======================================================
\noindent {\bf Abstract.} Recently, George
  Beck introduced
  two partition  statistics
   $NT(m,j,n)$ and $M_{\omega}(m,j,n)$,
    which denote
     the total number of parts in the partition
      of $n$ with rank congruent
       to $m$ modulo $j$ and
    the total number of ones
      in
    the partition
      of $n$ with crank congruent
       to $m$ modulo $j$, respectively.
        Andrews proved a congruence on $NT(m,5,n)$
         which was conjectured by Beck.
         Very recently, Chan, Mao and Osburn
           established
           a number of Andrews-Beck type
            congruences and posed
             several conjectures
              involving $NT(m,j,n)$
             and $M_{\omega}(m,j,n)$.
              Some of those  conjectures
              were proved by Chern and Mao.
               In this paper, we confirm
                the remainder three
                conjectures of Chan-Mao-Osburn
                 and two conjectures
                  due to Mao.
                 We also present two
                  new conjectures on
                   $M_{\omega}(m,j,n)$
                    and $NT(m,j,n)$.



\section{Introduction}

\allowdisplaybreaks

A partition
 $\pi=(\pi_1, \pi_2,\ldots,\pi_k)$
  of a positive integer
 $n$ is a   sequence
  of positive integers
   such that $\pi_1\geq \pi_2\geq \cdots
  \geq \pi_k>0$ and  $\pi_1+\pi_2+\cdots
  +\pi_k=n$. The $\pi_i$
   are called the parts of the partition
    \cite{Andrews-1976}.
    We shall write $\pi\vdash
     n$ to denote $\pi$ is a partition
     of $n$. We also use $\#(\pi)$ and
     $\lambda(\pi)$
      to denote the number of
       parts of $\pi$ and the
        largest part of $\pi$, respectively.
 As usual,
   let $p(n)$ denote the number of
    partitions of $n$ and set $p(0)=1$.
  In the theory of partition,
   one of the most well-known results
    is achieved by
      Ramanujan \cite{Ramanujan-1}. In
       1919, he
 found  that for $n\geq 0$,
 \begin{align*}
     p(5n + 4) &\equiv  0 \pmod 5,   \\
     p(7n + 5) & \equiv 0
      \pmod 7,    \\
      p(11n + 6) & \equiv 0
       \pmod {11}.
    \end{align*}

In 1944,   to  give
 combinatorial interpretations of Ramanujan's
 congruences, Dyson \cite{Dyson}
 defined   the rank of a partition
 to be  the largest part of the partition
minus the number of parts, namely,
\[
{\rm rank}(\pi):=\lambda(\pi)-\#(\pi).
\]
 For example, the rank
 of the partition $4+2+2+2
  +1$ is $4-5=-1$.
Let $N(m,j,n)$ count the number
 of partitions of $n$ with rank congruent
  to $m$ modulo $j$.
 Dyson \cite{Dyson} also conjectured
 that for $0\leq m \leq 4$
\begin{align}\label{1-1}
 N(m,5,5n+4)= \frac{ p(5n + 4)}{5},
\end{align}
and for $0\leq m \leq 6$,
\begin{align}\label{1-2}
 N(m,7,7n+5)= \frac{ p(7n + 5)}{7}
  .
\end{align}

 In 1954,   Atkin and
  Swinnerton-Dyer
  \cite{Atkin} proved
Dyson's conjectures \eqref{1-1}
 and \eqref{1-2}.
   Therefore Dyson's rank
gives combinatorial interpretations for Ramanujan's first two
congruences. Unfortunately,
  it turned out that Dyson's rank fails
 to  explain Ramanujan's
   third congruence modulo 11
    combinatorially. So Dyson conjectured
     the existence of an unknown
      partition statistic,
which he whimsically called ``the crank"
  to explain the third congruences modulo 11.
 In 1988, Andrews and Garvan
\cite{Andrews-2} finally
  found the actual crank.
For a partition
  $\pi$, let
   $\omega(\pi)$ denote the number of ones,
    and $\mu(\pi)$ the number
    of parts larger than $\omega(\pi)$. Then,
     the crank of $\pi$ is defined as follows
     \begin{align*}
{\rm crank}({\pi}):= \left\{
  \begin{aligned} &  \lambda(\pi)
   ,\qquad \qquad \quad {\rm if\ }
   \omega(\pi)=0,
   \\[6pt]
         & \mu(\pi)-\omega(\pi),\  \qquad {\rm otherwise }.
  \end{aligned} \right.
\end{align*}

 Recently, Andrews
 \cite{Andrews} mentioned that George
  Beck introduced
  two partition  statistics
   $NT(m,j,n)$ and $M_{\omega}(m,j,n)$,
    which denote
     the total number of parts in the partition
      of $n$ with rank congruent
       to $m$ modulo $j$ and
    the total number of ones
      in
    the partition
      of $n$ with crank congruent
       to $m$ modulo $j$, respectively, i.e.,
\[
NT(m,j,n)=\sum_{\pi \vdash n, \atop
 {\rm rank}(\pi)\equiv m
 \ ({\rm mod}\ j)}\#(\pi)
\]
and
\[
M_{\omega}(m,j,n)=\sum_{\pi
\vdash n, \atop
 {\rm crank}(\pi)\equiv m \
 ({\rm mod}\ j)}\omega(\pi).
\]
 Andrews \cite{Andrews}
  proved the following  Andrews-Beck
 type
  congruence which  was
   conjectured by Beck
\[
\sum_{m=1}^4
 m NT(m,5,5n+1) \equiv \sum_{m=1}^4
 m NT(m,5,5n+4) \equiv
0 \pmod 5.
\]

Motivated by Andrews's work, Chern
\cite{Chern-1,Chern-2,Chern-3}
proved some
 identities involving the weighted
  rank and crank moments and established
   a number of new Andrews-Beck type
   congruences on $NT(m,j,n)$
    and $M_{\omega}(m,j,n)$.
   For example, Chern \cite{Chern-1} proved that for
    $n\geq 0$,
\begin{align}\label{1-3}
\sum_{m=1}^4 mM_{\omega}(m,5,5n+4)\equiv
 0 \pmod 5.
\end{align}
In a recent paper,
 Lin, Peng and Toh
 \cite{Lin} considered the generalized
 crank defined
  by Fu and Tang \cite{Fu}
   for $k$-colored partitions
    and derived a number of Andrews-Beck
     type congruences. Very recently,
       Chan, Mao and Osburn
    \cite{Chan}
     proved
       three variations of Andrews-Beck type
        congruences and posed  a number of
         conjectures on $NT(m,j,n)$
          and $M_{\omega}(m,j,n)$.
 Those conjectures on Andrews-Beck type
  congruences of $NT(m,j,n)$
          and $M_{\omega}(m,j,n)$ were proved
           by Chern \cite{Chern-3}. Later,
            Mao \cite{mao} proved the following two identities
             on $NT(m,j,n)$ which were
              conjectured by Chan,
               Mao and Osburn \cite{Chan}:
\begin{align*}
&\sum_{n=0}^\infty
 \left(NT(1,7,7n+5)-NT(6,7,7n+5)
 +3NT(2,7,7n+5)-3NT(5,7,7n+5)\right)q^n
 \nonumber\\[6pt]
 =&-7\frac{(q^3,q^4,q^7,q^7,q^7;q^7)_\infty
  }{(q,q^2,q^2,q^5,q^5,q^6;q^7)_\infty}
\end{align*}
and
\begin{align*}
&\sum_{n=0}^\infty
 \left(NT(1,7,7n+4)-NT(6,7,7n+4)
 +2NT(3,7,7n+4)-2NT(4,7,7n+4)\right)q^n
 \nonumber\\[6pt]
 =&-7\frac{(q^3,q^3,q^4,q^4, q^7,q^7,q^7;q^7)_\infty
  }{(q,q^2,q^2,q^2,q^5,q^5,q^5,
  q^6;q^7)_\infty},
\end{align*}
 where here and throughout
 the rest of the paper, we adopt the standard
  $q$-series notation
  \[
(a;q)_\infty=\prod_{n=0}^\infty (1-aq^n)
  \]
and for each positive integer $k$,
\[
(a_1,a_2,\ldots,a_k;q)_\infty
 =(a_1;q)_\infty(a_2;q)_\infty
  \cdots (a_k;q)_\infty.
\]
In another paper, Mao \cite{mao-1}
 also gave several conjectures on $NT(m,j,n)$
  and $M_{\omega}(m,j,n)$.

The aim of the paper is to confirm  the
 remainder
three conjectures
 of Chan-Mao-Osburn \cite{Chan}
  and two conjectures due to Mao \cite{mao-1}
  on some relations
  involving
    $NT(m,5,n)$ and $M_{\omega}(m,5,n)$.

\begin{theorem}\label{Th-1}
We have
\begin{align}\label{1-4}
\sum_{n=0}^\infty&
 (NT(1,5,5n+4)-NT(4,5,5n+4)\nonumber\\[6pt]
 &\qquad +2M_{\omega}
 (2,5,5n+4)-2M_{\omega}(3,5,5n+4))q^n
 =-5\frac{(q^5;q^5)_\infty^4
 }{(q;q)_\infty}
\end{align}
and for $n\geq 0$,
\begin{align}\label{1-4-1}
M_{\omega}(2,5,5n+4)
-M_{\omega}(3,5,5n+4)  =
  2 NT(1,5,5n+4)- 2
NT(4 ,5,5n+4).
\end{align}
\end{theorem}

\begin{theorem}\label{Th-2}
 For $n\geq 0$,
\begin{align}\label{1-5}
M_{\omega}(1,5,5n+4)-M_{\omega}(4,5,5n+4)  =
  2 M_{\omega}(3,5,5n+4)- 2
 M_{\omega}(2 ,5,5n+4).
\end{align}
\end{theorem}

\begin{theorem}\label{Th-3}
For $n\geq 0$,
\begin{align}\label{1-6}
M_{\omega}(1,5,5n+2)
 -M_{\omega}(4,5,5n+2)=2
  NT(3,5,5n+2)-2NT(2,5,5n+2).
\end{align}
\end{theorem}
\begin{theorem}\label{Th-4}
For $n\geq 0$,
\begin{align}\label{1-6-1}
M_{\omega}(2,5,5n+1)
 -M_{\omega}(3,5,5n+1)=
  NT(2,5,5n+1)- NT(3,5,5n+1).
\end{align}
\end{theorem}

 Identities \eqref{1-4},
  \eqref{1-5}
   and \eqref{1-6} were
 first
   conjectured by Chan, Mao and Osburn
   \cite{Chan}
    and \eqref{1-4-1} and \eqref{1-6-1}
     were conjectured by Mao \cite{mao-1}.
     Moreover,
   identity \eqref{1-5}
   implies \eqref{1-3}.

\section{Preliminaries}

In this section, we present several lemmas
 which will be used
  to prove the main results of this paper.

The following lemma was given by Garvan \cite{Garvan}.

\begin{lemma} \label{L-1}
\cite[(3.1)]{Garvan}
 Let $\zeta={e^{2\pi i/5}}$. For $m=1,2$,
\begin{align}\label{2-1}
\frac{(q;q)_\infty
  }{(\zeta^m q;q)_\infty
  (q/\zeta^m;q)_\infty}
  &=A(q^5)-(\zeta^m+\zeta^{-m})^2
   qB(q^5)\nonumber\\[6pt]
   &\qquad +(\zeta^{2m}+\zeta^{-2m})
   q^2C(q^5)-(\zeta^m+\zeta^{-m})
    q^3D(q^5),
\end{align}
where
\begin{align}\label{2-2}
A(q)=&\frac{(q^2,q^3,q^5;q^5)_\infty
 }{(q,q^4;q^5)_\infty^2}, \
 B(q)= \frac{( q^5;q^5)_\infty
 }{(q,q^4;q^5)_\infty},\
 C(q)=\frac{( q^5;q^5)_\infty
 }{(q^2,q^3;q^5)_\infty}, \
 D(q)=\frac{(q,q^4,q^5;q^5)_\infty
 }{(q^2,q^3;q^5)_\infty^2}.
\end{align}

\end{lemma}

\begin{lemma} \label{L-2}
We have
\begin{align}
\sum_{n=0}^\infty
 \frac{q^n}{1-q^{5n+1}}-
 \sum_{n=0}^\infty
 \frac{q^{4n+3}}{1-q^{5n+4}}
&=\frac{(q^2,q^3,q^5,q^5;q^5)_\infty
  }{(q,q,q^4,q^4;q^5)_\infty}
,\label{2-3}\\[6pt]
  \sum_{n=0}^\infty
 \frac{q^{2n+1}}{1-q^{5n+3}}-
  \sum_{n=0}^\infty
 \frac{q^{3n+1}}{1-q^{5n+2}}
 &
 =0,\label{2-4}\\[6pt]
\sum_{n=0}^\infty
 \frac{q^n}{1-q^{5n+2}}
- \sum_{n=0}^\infty
 \frac{q^{4n+2}}{1-q^{5n+3}}& = \frac{(q^5,q^5;q^5)_\infty
  }{(q,q^4;q^5)_\infty}
,\label{2-5}\\[6pt]
 \sum_{n=0}^\infty
 \frac{q^{2n }}{1-q^{5n+2}}-
   \sum_{n=0}^\infty
 \frac{q^{3n+1}}{1-q^{5n+3}}
&
 =\frac{(q,q^4,q^5,q^5;q^5)_\infty
  }{(q^2,q^2,q^3,q^3;q^5)_\infty}
, \label{2-6}
\\[6pt]
 \sum_{n=0}^\infty
 \frac{q^n}{1-q^{5n+3}}
- \sum_{n=0}^\infty
 \frac{q^{4n+1}}{1-q^{5n+2}}&=\frac{( q^5,q^5;q^5)_\infty
  }{( q^2,q^3 ;q^5)_\infty}
,\label{2-7}\\[6pt]
 \sum_{n=0}^\infty
 \frac{q^{2n+1 }}{1-q^{5n+4}}-
  \sum_{n=0}^\infty
 \frac{q^{3n}}{1-q^{5n+1}}
&=-\frac{( q^5,q^5;q^5)_\infty
  }{( q^2,q^3 ;q^5)_\infty}
,\label{2-8}\\[6pt]
\sum_{n=0}^\infty
 \frac{q^n}{1-q^{5n+4}}
- \sum_{n=0}^\infty
 \frac{q^{4n}}{1-q^{5n+1}}&=0
,\label{2-9}\\[6pt]
\sum_{n=0}^\infty
 \frac{q^{2n }}{1-q^{5n+1}}-
  \sum_{n=0}^\infty
 \frac{q^{3n+2}}{1-q^{5n+4}}&
 =\frac{(q^5,q^5;q^5)_\infty
  }{(q,q^4;q^5)_\infty}
,\label{2-10}\\[6pt]
 \sum_{n=0}^\infty
 \frac{q^{2n+1 }}{1-q^{5n+2}}-
  \sum_{n=0}^\infty
 \frac{q^{3n+2}}{1-q^{5n+3}}
&=q\frac{(q,q^4,q^5,q^5;q^5)_\infty
  }{(q^2,q^2,q^3,q^3;q^5)_\infty}
.\label{v-1}
  \end{align}
\end{lemma}

\noindent{\it Proof.} Let
 $r,s,t$ be integers with $1\leq r,s\leq 4$
 and $0\leq t\leq 4$.  It is easy to check
 that
 \begin{align}
\sum_{n=0}^\infty \frac{q^{rn+t}} {1-q^{5n+s}} -\sum_{n=0}^\infty
 \frac{q^{(5-r)n+5+t-r-s}}{1-q^{5n+5-s}}
 =&\sum_{n=0}^\infty \frac{q^{rn+t}}
 {1-q^{5n+s}} -\sum_{n=1}^\infty
 \frac{q^{(5-r)(n-1)+5+t-r-s}}{1-q^{5n -s}}
\nonumber\\[6pt]
=&\sum_{n=0}^\infty \frac{q^{rn+t}}
 {1-q^{5n+s}} -\sum_{n=-\infty}^{-1}
 \frac{q^{(5-r)(-n-1)+5+t-r-s}}{1-q^{-5n -s}}
\nonumber\\[6pt]
=&\sum_{n=0}^\infty \frac{q^{rn+t}}
 {1-q^{5n+s}} +\sum_{n=-\infty}^{-1}
 \frac{q^{rn+t}}{1-q^{ 5n
  +s}}
\nonumber\\[6pt]
=&\sum_{n=-\infty}^\infty \frac{q^{rn+t}}
 {1-q^{5n+s}}.\label{2-11}
 \end{align}
It follows from \cite[Lemma 4.4.2, p. 117]{Andrews-0}
  that
\begin{align}\label{2-12}
\sum_{n=-\infty}^\infty
 \frac{q^{ni}}{1-q^{5n+j}}
 =\frac{(q^{i+j},q^{5-i-j},q^5,q^5;q^5)_\infty
  }{(q^i,q^j,q^{5-i},q^{5-j};q^5)_\infty}.
\end{align}
Combining \eqref{2-11}
 and \eqref{2-12} yields
\begin{align}\label{2-13}
\sum_{n=0}^\infty \frac{q^{rn+t}} {1-q^{5n+s}} -\sum_{n=0}^\infty
 \frac{q^{(5-r)n+5+t-r-s}}{1-q^{5n+5-s}}
 =q^t\frac{(q^{r+s},q^{5-r-s},q^5,q^5;q^5)_\infty
  }{(q^r,q^s,q^{5-r},q^{5-s};q^5)_\infty}.
\end{align}
The proofs below make frequent use of
 \eqref{2-13}.

To prove \eqref{2-3}, use \eqref{2-13}
 with $r=1,s=1$ and $t=0$.

 To prove \eqref{2-4}, use \eqref{2-13}
 with $r=2,s=3$ and $t=1$.

 To prove \eqref{2-5}, use \eqref{2-13}
 with $r=1,s=2$ and $t=0$.

  To prove \eqref{2-6}, use \eqref{2-13}
 with $r=2,s=2$ and $t=0$.

  To prove \eqref{2-7}, use \eqref{2-13}
 with $r=1,s=3$ and $t=0$.

   To prove \eqref{2-8}, use \eqref{2-13}
 with $r=2,s=4$ and $t=1$.

   To prove \eqref{2-9}, use \eqref{2-13}
 with $r=1,s=4$ and $t=0$.

   To prove \eqref{2-10}, use \eqref{2-13}
 with $r=2,s=1$ and $t=0$.

   To prove \eqref{v-1}, use \eqref{2-13}
 with $r=2,s=2$ and $t=1$.

The proof of Lemma \ref{L-2} is complete. \qed

\begin{lemma} \label{L-3}
We have
\begin{align}\label{2-21}
\sum_{n=1}^\infty
     \frac{q^n+q^{2n}-q^{3n}-q^{4n}}{1-q^{5n}}
= \frac{ 2(q^2,q^3,q^5;q^5)_\infty^2
 }{5(q,q^4;q^5)_\infty^3}-  \frac{q
 (q,q^4,q^5;q^5)_\infty^2
  }{5(q^2,q^3;q^5)_\infty^3}-\frac{2}{5}
\end{align}
and
\begin{align}\label{2-21-1}
\sum_{n=1}^\infty
     \frac{q^n-2q^{2n}+2q^{3n}-q^{4n}}{1-q^{5n}}
= \frac{ (q^2,q^3,q^5;q^5)_\infty^2
 }{10(q,q^4;q^5)_\infty^3}+  \frac{7q
 (q,q^4,q^5;q^5)_\infty^2
  }{10(q^2,q^3;q^5)_\infty^3}-\frac{1}{10}.
\end{align}
\end{lemma}

\noindent{\it Proof.}
  It is easy to verify that for $1\leq m \leq
4$,
\begin{align}\label{2-22}
\sum_{n=1}^\infty \frac{q^{mn}}{1-q^{5n}}
 =&\sum_{n=1}^\infty q^{mn}
   \sum_{j=0}^\infty q^{5nj}
   \nonumber\\[6pt]
 =&\sum_{j=0}^\infty\sum_{n=1}^\infty
  q^{(5j+m)n}= \sum_{j=0}^\infty
\frac{q^{5j+m}}{1-q^{5j+m}}.
\end{align}
Thanks to \eqref{2-22},
\begin{align}\label{2-23}
\sum_{n=1}^\infty
     \frac{q^n+q^{2n}-q^{3n}-q^{4n}}{1-q^{5n}}
=\sum_{j=0}^\infty \left(
 \frac{q^{5j+1}}{1-q^{5j+1}}
-\frac{q^{5j+4}}{1-q^{5j+4}} +\frac{q^{5j+2}}{1-q^{5j+2}}
-\frac{q^{5j+3}}{1-q^{5j+3}}\right)  .
\end{align}
It is easy to check that
\begin{align}\label{2-24}
\sum_{j=0}^\infty \left(
 \frac{q^{5j+1}}{1-q^{5j+1}}
-\frac{q^{5j+4}}{1-q^{5j+4}}\right)= \sum_{j=0}^\infty \left(
 \frac{q^{5j+6}}{1-q^{5j+6}}
-\frac{q^{5j-1}}{1-q^{5j-1}}\right)-1
\end{align}
and
\begin{align}\label{2-25}
\sum_{j=0}^\infty \left(
 \frac{q^{5j+2}}{1-q^{5j+2}}
-\frac{q^{5j+3}}{1-q^{5j+3}}\right)=
 \sum_{j=0}^\infty \left(
 \frac{q^{5j+7}}{1-q^{5j+7}}
-\frac{q^{5j-2}}{1-q^{5j-2}}\right)-1.
\end{align}
It follows from \cite[Lemma 2.4, (2.11) and (2.12)]{mao} that
\begin{align}\label{2-26}
 \sum_{j=0}^\infty \left(
 \frac{q^{5j+6}}{1-q^{5j+6}}
-\frac{q^{5j-1}}{1-q^{5j-1}}\right) =\frac{1}{10} \left(
\frac{3(q^2,q^3,q^5;q^5)_\infty^2
 }{(q,q^4;q^5)_\infty^3}+q\frac{
 (q,q^4,q^5;q^5)_\infty^2
  }{(q^2,q^3;q^5)_\infty^3}+7\right)
\end{align}
and
\begin{align}\label{2-27}
 \sum_{j=0}^\infty \left(
 \frac{q^{5j+7}}{1-q^{5j+7}}
-\frac{q^{5j-2}}{1-q^{5j-2}}\right)= \frac{1}{10} \left( \frac{
(q^2,q^3,q^5;q^5)_\infty^2
 }{(q,q^4;q^5)_\infty^3}-3q\frac{
 (q,q^4,q^5;q^5)_\infty^2
  }{(q^2,q^3;q^5)_\infty^3}+9\right)
\end{align}
Combining \eqref{2-23}--\eqref{2-27}, we arrive at
 \eqref{2-21}. This completes the proof of Lemma
 \ref{L-3}.  \qed

\section{The generating functions for
 $M_{\omega}(a,5,n)$}

 In this section, we establish
  the generating functions for
   $M_{\omega}(a,5,n)$.

 \begin{theorem} \label{Th-3-1}
  We have
  \begin{align}
\sum_{n\geq 0} M_{\omega}(0,5,n)q^n
 =&\frac{q^3}{5} D(q^5)
\left(-3R_1(q)+2R_2(q)+2R_3(q) -3R_4(q)+2 R_5(q)-2S(q)\right)
\nonumber\\[6pt]
&+\frac{q^2}{5} C(q^5) \left(-2R_1(q)+3 R_2(q)+3R_3(q)
-2R_4(q)-2R_5(q)+2S(q)\right)
\nonumber\\[6pt]
&+\frac{q}{5} B(q^5) \left(4R_1(q)-
  R_2(q)-R_3(q)
+4R_4(q)-6R_5(q)+6S(q)\right)
\nonumber\\[6pt]
&+\frac{1}{5} A(q^5) \left(-R_1(q)-
  R_2(q)-R_3(q)
-R_4(q)+4R_5(q)-4S(q)\right)+T(q),
\label{3-1}\\[6pt]
\sum_{n\geq 0} M_{\omega}(1,5,n)q^n
 =&\frac{q^3}{5} D(q^5)
\left(2R_1(q)+2R_2(q)-3R_3(q) +2R_4(q)-3
  R_5(q)+3S(q)\right)
\nonumber\\[6pt]
&+\frac{q^2}{5} C(q^5) \left(3
 R_1(q)+3 R_2(q)-2R_3(q)
-2R_4(q)-2R_5(q)+2S(q)\right)
\nonumber\\[6pt]
&+\frac{q}{5} B(q^5) \left(-R_1(q)-
  R_2(q)+4R_3(q)
-6R_4(q)+4R_5(q)-4S(q)\right)
\nonumber\\[6pt]
&+\frac{1}{5} A(q^5) \left(-R_1(q)-
  R_2(q)-R_3(q)
+4R_4(q)-R_5(q)+S(q)\right)+T(q),
\label{3-2}\\[6pt]
\sum_{n\geq 0} M_{\omega}(2,5,n)q^n
 =&\frac{q^3}{5} D(q^5)
\left(2R_1(q)-3R_2(q)+2R_3(q) -3R_4(q)+2R_5(q)-2S(q)\right)
\nonumber\\[6pt]
&+\frac{q^2}{5} C(q^5) \left(3
 R_1(q)-2 R_2(q)-2R_3(q)
-2R_4(q)+3R_5(q)-3S(q)\right)
\nonumber\\[6pt]
&+\frac{q}{5} B(q^5) \left(-R_1(q)+4 R_2(q)-6R_3(q)
+4R_4(q)-R_5(q)+S(q)\right)
\nonumber\\[6pt]
&+\frac{1}{5} A(q^5) \left(-R_1(q)-
  R_2(q)+4R_3(q)
-R_4(q)-R_5(q)+S(q)\right)+T(q),
\label{3-3}\\[6pt]
\sum_{n\geq 0} M_{\omega}(3,5,n)q^n
 =&\frac{q^3}{5} D(q^5)
\left(-3R_1(q)+2R_2(q)-3R_3(q) +2
 R_4(q)+2R_5(q)-2S(q)\right)
\nonumber\\[6pt]
&+\frac{q^2}{5} C(q^5) \left(-2
 R_1(q)-2 R_2(q)-2R_3(q)
+3R_4(q)+3R_5(q)-3S(q)\right)
\nonumber\\[6pt]
&+\frac{q}{5} B(q^5) \left(4R_1(q)
 -6 R_2(q)+4R_3(q)
-R_4(q)-R_5(q)+S(q)\right)
\nonumber\\[6pt]
&+\frac{1}{5} A(q^5) \left(-R_1(q)+4
  R_2(q)-R_3(q)
-R_4(q)-R_5(q)+S(q)\right)+T(q),
\label{3-4}\\[6pt]
\sum_{n\geq 0} M_{\omega}(4,5,n)q^n
 =&\frac{q^3}{5} D(q^5)
\left(2R_1(q)-3R_2(q)+2R_3(q) +2R_4(q)-3
 R_5(q)+3S(q)\right)
\nonumber\\[6pt]
&+\frac{q^2}{5} C(q^5) \left(-2
 R_1(q)-2 R_2(q)+3R_3(q)
+3R_4(q)-2R_5(q)+2S(q)\right)
\nonumber\\[6pt]
&+\frac{q}{5} B(q^5) \left(-6
 R_1(q)+4 R_2(q)-R_3(q)
-R_4(q)+4R_5(q)-4S(q)\right)
\nonumber\\[6pt]
&+\frac{1}{5} A(q^5) \left(4R_1(q)- R_2(q)-R_3(q)
-R_4(q)-R_5(q)+S(q)\right)+T(q),\label{3-5}
\end{align}
where $A(q)$, $B(q)$, $C(q)$, $D(q)$
 are defined by \eqref{2-2} and
\begin{align}\label{3-6}
R_i(q)=\sum_{n=1}^\infty \frac{q^{ni} }{1-q^{5n}}, \qquad
S(q)=\sum_{n=1}^\infty \frac{q^{n+1}}{1-q^{n+1}}, \qquad T(q)
=\frac{q}{5(1-q)(q;q)_\infty}.
\end{align}

 \end{theorem}

 \noindent{\it Proof.}
Chern \cite[(3.2)]{Chern-1} proved that
\begin{align}\label{3-7}
\sum_{n\geq 0} \sum_{\lambda
 \vdash n}   \omega(\lambda)
z^{crank(\lambda)}  q^n=\frac{ ( q ;q)_\infty
    }{(zq,xq/z;q)_\infty}\sum_{n\geq 1}
    \left(\frac{q^n/z}{1-q^n/z}-\frac{q^{n+1}}{
    1-q^{n+1}}\right).
\end{align}
By \eqref{3-7} and the definition
 of $M_{\omega}(b,5,n)$,
 \begin{align}
\sum_{n\geq 0} M_{\omega}(b,5,n)q^n = &
\frac{1}{5} \sum_{j=0}^4
 \zeta^{-bj}\frac{(q;q)_\infty
  }{(\zeta^j q;q)_\infty
  (q/\zeta^j;q)_\infty}
   \left(\sum_{n=1}^\infty
    \frac{\zeta^{-j}q^n}{1-q^n\zeta^{-j}}
     -S(q)\right)\nonumber\\[6pt]
     =& T(q)
     +\frac{1}{5} \sum_{j=1}^4
 \zeta^{-bj}\frac{(q;q)_\infty
  }{(\zeta^j q;q)_\infty
  (q/\zeta^j;q)_\infty}
   \left(\sum_{n=1}^\infty
    \frac{\zeta^{-j}q^n}{1-q^n\zeta^{-j}}
     -S(q)\right),\label{3-8}
 \end{align}
where  $S(q)$ and $T(q)$ are
 defined by \eqref{3-6} and $\zeta=e^{2\pi i/5}$.
 Moreover, it is easy to check that
\begin{align}\label{3-9}
\sum_{n=1}^\infty
 \frac{\zeta^{-j} q^n}{1-\zeta^{-j}q^{n}}
  = &\sum_{n=1}^\infty
   \frac{q^{5n}}{1-q^{5n}}+\zeta^{-j}
   \sum_{n=1}^\infty
   \frac{q^{ n}}{1-q^{5n}}
   +\zeta^{-2j}
   \sum_{n=1}^\infty
   \frac{q^{ 2n}}{1-q^{5n}}
   \nonumber\\[6pt]
&    +\zeta^{-3j}
   \sum_{n=1}^\infty
   \frac{q^{ 3n}}{1-q^{5n}}
    +\zeta^{-4j}
   \sum_{n=1}^\infty
   \frac{q^{ 4n}}{1-q^{5n}}.
\end{align}
Setting $b=0,1,2,3,4$ in \eqref{3-8}
 and substituting \eqref{2-1} and
  \eqref{3-9}
   into \eqref{3-8},
  we arrive at \eqref{3-1}--\eqref{3-5}, respectively.
   This completes the proof of Theorem
   \ref{Th-3-1}. \qed

\section{Proofs of Theorems
\ref{Th-1}--\ref{Th-4}}

The objective of this section is to prove
 Theorems \ref{Th-1}--\ref{Th-4}.

\noindent{\it Proof of Theorem \ref{Th-1}.} In light of \eqref{3-3}
and \eqref{3-4},
\begin{align}\label{4-1}
&\sum_{n\geq 0}(   M_{\omega}(2 ,5,n)-
 M_{\omega}(3,5,n)
)q^n \nonumber\\[6pt]
=& q^3D(q^5)(R_1(q)-R_2(q)+R_3(q)-R_4(q))
+q^2C(q^5)(R_1(q)-R_4(q)) \nonumber\\[6pt]
&+qB(q^5)(-R_1(q)+2R_2(q)-2R_3(q)+R_4(q)) -A(q^5)( R_2(q)- R_3(q) )
\end{align}
If we extract
 those  terms in which the power of $q$
 is
  congruent to 4 modulo 5 in \eqref{4-1}, then
   divided by $q^4$
  and
   replace  $q^5$ by $q$, we arrive at
\begin{align}\label{4-2}
&\sum_{n\geq 0}(  M_{\omega}(2 ,5,5n+4)-  M_{\omega}(3,5,5n+4)
  )q^n \nonumber\\[6pt]
=&D(q) \left(  \sum_{n=0}^\infty
 \frac{q^n}{1-q^{5n+1}}-
  \sum_{n=0}^\infty
 \frac{q^{4n+3}}{1-q^{5n+4}}\right)
 -D(q)\left(\sum_{n=0}^\infty
 \frac{q^{2n+1}}{1-q^{5n+3}}-
  \sum_{n=0}^\infty
 \frac{q^{3n+1}}{1-q^{5n+2}}\right)
 \nonumber\\[6pt]
& +C(q) \left(  \sum_{n=0}^\infty
 \frac{q^n}{1-q^{5n+2}}
- \sum_{n=0}^\infty
 \frac{q^{4n+2}}{1-q^{5n+3}}
  \right)-B(q) \left(  \sum_{n=0}^\infty
 \frac{q^n}{1-q^{5n+3}}
- \sum_{n=0}^\infty
 \frac{q^{4n+1}}{1-q^{5n+2}}\right)
  \nonumber\\[6pt]
&
 +2B(q)\left(\sum_{n=0}^\infty
 \frac{q^{2n+1 }}{1-q^{5n+4}}-
  \sum_{n=0}^\infty
 \frac{q^{3n}}{1-q^{5n+1}}\right)
  -A(q) \left( \sum_{n=0}^\infty
 \frac{q^{2n }}{1-q^{5n+2}}-
   \sum_{n=0}^\infty
 \frac{q^{3n+1}}{1-q^{5n+3}}\right).
\end{align}
Thanks to \eqref{2-3}--\eqref{2-8} and \eqref{4-2},
\begin{align}\label{4-3}
 \sum_{n=0}^\infty (M_{\omega} (2,5,5n+4) -
  M_{\omega} (3,5,5n+4)) q^n&=
  -2\frac{(q^5;q^5)_\infty^4}{(q;q)_\infty}
\end{align}
Very recently, Mao \cite[Theorem 1.1]{mao} proved the  following
identity:
\begin{align}\label{4-4}
\sum_{n=0}^\infty (NT (1,5,5n+4) -
  NT(4,5,5n+4)) q^n&=
  -\frac{(q^5;q^5)_\infty^4}{(q;q)_\infty}.
\end{align}
Theorem \ref{Th-1}
 follows from \eqref{4-3} and \eqref{4-4}.

Now, we turn to prove Theorem \ref{Th-2}.

\noindent{\it Proof of Theorem \ref{Th-2}.}
      In view of \eqref{3-2} and \eqref{3-5},
 \begin{align}
 &\sum_{n=0}^\infty (M_{\omega} (1,5,n) -
  M_{\omega} (4,5,n)) q^n
  \nonumber\\[6pt]
  =&q^3D(q^5)(R_2(q)-R_3(q))+q^2C(q^5)
  (R_1(q)-R_4(q)+R_2(q)-R_3(q))\nonumber\\[6pt]
  &+qB(q^5)(R_1(q)-R_4(q)-R_2(q)+R_3(q))
   -A(q^5)(R_1(q)-R_4(q)).\label{4-5}
\end{align}
Extracting
 those  terms in which the power of $q$ is
  congruent to 4 modulo 5 in \eqref{4-5}, then
  dividing by $q^4$
  and
   replacing  $q^5$ by $q$, we deduce that
\begin{align}\label{4-6}
 &\sum_{n=0}^\infty ( M_{\omega} (1,5,5n+4) -
  M_{\omega} (4,5,5n+4)) q^n
  \nonumber\\[6pt]
  =&D(q) \left(
  \sum_{n=0}^\infty
 \frac{q^{2n+1}}{1-q^{5n+3}}-
  \sum_{n=0}^\infty
 \frac{q^{3n+1}}{1-q^{5n+2}}\right)
 +C(q) \left(  \sum_{n=0}^\infty
 \frac{q^n}{1-q^{5n+2}}
-  \sum_{n=0}^\infty
 \frac{q^{4n+2}}{1-q^{5n+3}}\right)
 \nonumber\\[6pt]
&
 +C(q)\left(\sum_{n=0}^\infty
 \frac{q^{2n }}{1-q^{5n+1}}-
  \sum_{n=0}^\infty
 \frac{q^{3n+2}}{1-q^{5n+4}}\right)
  +B(q) \left(  \sum_{n=0}^\infty
 \frac{q^n}{1-q^{5n+3}}
- \sum_{n=0}^\infty
 \frac{q^{4n+1}}{1-q^{5n+2}}\right)
  \nonumber\\[6pt]
& -B(q)\left( \sum_{n=0}^\infty
 \frac{q^{2n+1 }}{1-q^{5n+4}}-
   \sum_{n=0}^\infty
 \frac{q^{3n}}{1-q^{5n+1}}\right)
  -A(q) \left(  \sum_{n=0}^\infty
 \frac{q^n}{1-q^{5n+4}}
- \sum_{n=0}^\infty
 \frac{q^{4n}}{1-q^{5n+1}}
 \right).
\end{align}
It follows from \eqref{2-4}, \eqref{2-5},
 \eqref{2-7}--\eqref{2-10}
  and \eqref{4-6} that
\begin{align}\label{4-7}
 \sum_{n=0}^\infty ( M_{\omega} (1,5,5n+4) -
  M_{\omega} (4,5,5n+4)) q^n&=
  4\frac{(q^5;q^5)_\infty^4}{(q;q)_\infty}.
\end{align}
Identity \eqref{1-5} follows from \eqref{4-3}
  and \eqref{4-7} and the proof
 of Theorem \ref{Th-2} is complete.

Next, we present a proof
 of Theorem \ref{Th-3}.

 \noindent{\it Proof of Theorem \ref{Th-3}.}
If we pick out
 those  terms in which the power of
  $q$ is
  congruent to 2 modulo 5 in \eqref{4-5}, then
   divided by $q^2$ and
   replace $q^5$ by $q$, we obtain
\begin{align}\label{4-8}
 &\sum_{n=0}^\infty ( M_{\omega} (1,5,5n+2) -
  M_{\omega} (4,5,5n+2) ) q^n
  \nonumber\\[6pt]
  =&D(q) \left(
  \sum_{n=0}^\infty
 \frac{q^{2n+1}}{1-q^{5n+2}}-
  \sum_{n=0}^\infty
 \frac{q^{3n+2}}{1-q^{5n+3}}\right)
 +C(q) {\sum_{n=1}^\infty}
     \frac{q^n+q^{2n}-q^{3n}
     -q^{4n}}{1-q^{5n}}
     \nonumber\\[6pt]
     &+B(q) \left(  \sum_{n=0}^\infty
 \frac{q^n}{1-q^{5n+1}}
- \sum_{n=0}^\infty
 \frac{q^{4n+3}}{1-q^{5n+4}}\right)
 -B(q)\left( \sum_{n=0}^\infty
 \frac{q^{2n+1 }}{1-q^{5n+3}}-
   \sum_{n=0}^\infty
 \frac{q^{3n+1}}{1-q^{5n+2}}\right)
 \nonumber\\[6pt]
 &-A(q) \left(  \sum_{n=0}^\infty
 \frac{q^n}{1-q^{5n+2}}
- \sum_{n=0}^\infty
 \frac{q^{4n+2}}{1-q^{5n+3}}
 \right).
\end{align}
It follows from \eqref{2-3}--\eqref{2-5},
 \eqref{v-1},
  \eqref{2-21} and \eqref{4-8} that
\begin{align}\label{4-9}
\sum_{n=0}^\infty ( M_{\omega} (1,5,5n+2) -
  M_{\omega} (4,5,5n+2)) q^n&=
  \frac{4q}{5}\frac{(q,q^4;q^5)_\infty^2(q^5
  ;q^5)_\infty^3
   }{(q^2,q^3;q^5)_\infty^4}\nonumber\\[6pt]
   &\quad +\frac{2}{5}\frac{(q^2,q^3;q^5)_\infty
  (q^5;q^5)_\infty^3
  }{(q,q^4;q^5)_\infty^3}-\frac{2}{5}
  \frac{(q^5;q^5)_\infty}{
  (q^2,q^3;q^5)_\infty}
.
\end{align}
Mao \cite[Theorem 1.2]{mao} also proved that
\begin{align}\label{4-10}
2\sum_{n=0}^\infty
 (NT(2,5,5n+2)-NT(3,5,5n+2))q^n
 =&-\frac{2}{5}\frac{(q^2,q^3;q^5)_\infty
  (q^5;q^5)_\infty^3
  }{(q,q^4;q^5)_\infty^3}
  \nonumber\\[6pt]
  &+\frac{2}{5}
  \frac{(q^5;q^5)_\infty}{
  (q^2,q^3;q^5)_\infty}
  -\frac{4q}{5}\frac{(q,q^4;q^5)_\infty^2
(q^5 ;q^5)_\infty^3
   }{(q^2,q^3;q^5)_\infty^4},
\end{align}
which yields  \eqref{1-6} after combining \eqref{4-9}.
 The proof of Theorem \ref{Th-3} is
 complete.

 Finally, we give  a proof
 of Theorem \ref{Th-4}.

 \noindent{\it Proof of Theorem \ref{Th-4}.}
Extracting
 those  terms in which the power of $q$ is
  congruent to 1 modulo 5 in \eqref{4-1}, then
  dividing by $q^4$
  and
   replacing  $q^5$ by $q$, we arrive at
\begin{align}\label{4-11}
&\sum_{n=0}^\infty (
 M_{\omega}(2,5,5n+1)
  -M_{\omega}(3,5,5n+1)) q^n
  \nonumber\\[6pt]
  =&D(q)\left(\sum_{n\geq 0}
   \frac{q^{n+1}}{1-q^{5n+3}}
   -\sum_{n\geq 0}
   \frac{q^{4n+2}}{1-q^{5n+2}}\right)
   -D(q)\left(\sum_{n\geq 0}
   \frac{q^{2n+2}}{1-q^{5n+4}}
   -\sum_{n\geq 0}
   \frac{q^{3n+1}}{1-q^{5n+1}}\right)
   \nonumber\\[6pt]
   &+C(q)\left(\sum_{n\geq 0}
   \frac{q^{ n+1}}{1-q^{5n+4}}
   -\sum_{n\geq 0}
   \frac{q^{4n+1}}{1-q^{5n+1}}\right)
   -B(q)\sum_{n\geq 1} \frac{q^n-2q^{2n}+2q^{3n}-q^{4n}
    }{1-q^{5n}}
    \nonumber\\[6pt]
   &-A(q)\left(\sum_{n\geq 0}
   \frac{q^{ 2n+1}}{1-q^{5n+3}}
   -\sum_{n\geq 0}
   \frac{q^{3n+1}}{1-q^{5n+2}}\right).
\end{align}
Thanks to \eqref{2-2}, \eqref{2-4},
 \eqref{2-7}--\eqref{2-9},
  \eqref{2-21-1}
  and \eqref{4-11},
\begin{align}\label{4-12}
\sum_{n=0}^\infty (
 M_{\omega}(2,5,5n+1)
  -M_{\omega}(3,5,5n+1)) q^n=&
  \frac{(q^5;q^5)_\infty}{
 10 (q,q^4;q^5)_\infty}
  -\frac{(q^2,q^3;q^5)_\infty
   (q^5;q^5)_\infty^3}{10(q,q^4;q^5)_\infty^4
    }\nonumber\\[6pt]
    &+\frac{13q(q,q^4;q^5)_\infty
     (q^5;q^5)_\infty^3}{
     10(q^2,q^3;q^5)_\infty^3}.
\end{align}
Mao \cite[Theorem 1.1, (1.12)]{mao}
 proved that
\begin{align}\label{4-13}
\sum_{n=0}^\infty ( NT(2,5,5n+1)
  -NT(3,5,5n+1)) q^n=&
  \frac{(q^5;q^5)_\infty}{
 10 (q,q^4;q^5)_\infty}
  -\frac{(q^2,q^3;q^5)_\infty
   (q^5;q^5)_\infty^3}{10(q,q^4;q^5)_\infty^4
    }\nonumber\\[6pt]
    &+\frac{13q(q,q^4;q^5)_\infty
     (q^5;q^5)_\infty^3}{
     10(q^2,q^3;q^5)_\infty^3},
\end{align}
which yields \eqref{1-6-1}
 after combining \eqref{4-12}. This
 completes the proof of Theorem \ref{Th-4}.
  \qed

\section{Concluding Remarks}

As seen in Introduction, the
  two partition  statistics
   $NT(m,j,n)$ and $M_{\omega}(m,j,n)$
    introduced by Beck have received
 a lot of attention in recent years.
 In particular, Chan, Mao and Osburn
 \cite{Chan} posed several conjectures on
$NT(m,j,n)$ and $M_{\omega}(m,j,n)$
 and some of them have been proved by
  Chern \cite{Chern-3} and Mao \cite{mao}.
 In this paper, we confirm
  the
 remainder
three conjectures
 of Chan-Mao-Osburn  \cite{Chan}
  and two conjectures given by Mao \cite{mao-1}
  on the relations
  between   $NT(m,j,n)$
  and $M_{\omega}(m,j,n)$
    by proving some $q$-series identities.
 From the results of this paper, we can obtain
  some congruences. For example,  by
 \eqref{4-3}, \eqref{4-7} and \eqref{4-9},
  we see that for $n\geq 0$,
 \begin{align}
M_{\omega}(2,5,5n+4)&\equiv M_{\omega}(3,5,5n+4)
 \pmod 2,\label{5-1}\\[6pt]
 M_{\omega}(1,5,5n+2)&\equiv
  M_{\omega}(4,5,5n+2)
 \pmod 2,\label{5-2}\\[6pt]
 M_{\omega}(1,5,5n+4)&\equiv M_{\omega}(4,5,5n+4)
 \pmod 4. \label{5-3}
 \end{align}
 It is interesting to give combinatorial
  proofs of  \eqref{5-1}--\eqref{5-3}.

Computer evidence suggests
 that the following two conjectures
might hold.

\begin{conjecture}
 Let $0\leq i<j\leq 4$ be integers. If
  $i+j\neq 5$, then
  \begin{align}\label{5-4}
\lim_{n\rightarrow \infty} \frac{\#\{k |M_{\omega}(i,5,k) \equiv
M_{\omega}(j,5,k) \pmod 2,\ 1\leq k \leq n
 \}}{n}=\frac{1}{2}.
  \end{align}
Moreover,
  \begin{align}\label{5-5}
\lim_{n\rightarrow \infty} \frac{\#\{k |M_{\omega}(1,5,k) \equiv
M_{\omega}(4,5,k) \pmod 2,\ 1\leq k \leq n
 \}}{n}=\frac{3}{10}
\end{align}
  and
  \begin{align}\label{5-6}
\lim_{n\rightarrow \infty} \frac{\#\{k |M_{\omega}(2,5,k) \equiv
M_{\omega}(3,5,k) \pmod 2,\ 1\leq k \leq n
 \}}{n}=\frac{2}{5}.
 \end{align}

\end{conjecture}

\begin{conjecture}
 Let $0\leq i<j\leq 4$ be integers. Then
  \begin{align}\label{5-7}
\lim_{n\rightarrow \infty} \frac{\#\{k |NT(i,5,k) \equiv NT(j,5,k)
\pmod 2,\ 1\leq k \leq n
 \}}{n}=\frac{1}{2}.
  \end{align}

\end{conjecture}

 \vspace{0.5cm}

 \subsection*{Acknowledgements}
The authors
     cordially thank  Renrong Mao    for his  helpful
     comments.
 This work was supported by
 the National Science Foundation
 of China (grant No.~11971203) and
    the Nature Foundation for Distinguished Young
   Scientists of Jiangsu Province
    (No.~BK20180044).

\end{document}